\documentclass{article}

\usepackage{amsmath}
\usepackage{amsfonts}
\usepackage{amssymb}
\usepackage{amscd}
\usepackage{amsbsy}
\usepackage{amsthm}
\usepackage{graphicx}
\usepackage{psfrag}
\usepackage{color}

\usepackage{cleveref}
\usepackage{mathrsfs}

\allowdisplaybreaks[3]

\usepackage[twoside]{geometry}
\newfont{\myfnt}{cmssi10 scaled 1440}
\numberwithin{equation}{section}
\catcode`@=11
\def\ps@nk{\def\@oddhead{\vbox{\hbox to \hsize{\pic \footnotesize \it \shorttitle
				\hfill \rm \thepage} \vspace{1mm} \vspace*{-2mm}}}
	\def\@evenhead{\vbox{\hbox to \hsize{\pic \footnotesize \rm \thepage \hfill \it \shortauthor}
			\vspace{1mm} \vspace*{-2mm}}}
	\def\@oddfoot{} \def\@evenfoot{}}
\def\ps@first{\def\@oddhead{}
	\def\@evenhead{}
	\def\@oddfoot{} \def\@evenfoot{}}
\def\ps@total{\def\@oddhead{\vbox{\hbox to \hsize{\footnotesize \rm \hfill TOTAL\ \ CONTENTS
				\hfill \thepage} \vspace{1mm} \hrule \vspace*{-2mm}}}
	\def\@evenhead{\vbox{\hbox to \hsize{\footnotesize \rm \thepage \hfill CHIN.\ \ ANN.\ \ MATH.
				\hfill} \vspace{1mm} \hrule \vspace*{-2mm}}}
	\def\@oddfoot{} \def\@evenfoot{}}
\newtheoremstyle{thmstyle}% name
{6pt}%      Space above
{6pt}%      Space below
{\it}%         Body font
{}%         Indent amount (empty = no indent, \parindent = para indent)
{\bf}% Thm head font
{}%        Punctuation after thm head
{.5em}%     Space after thm head: " " = normal interword space;
{}%         Thm head spec (can be left empty, meaning `normal')
\newtheoremstyle{remstyle}% name
{6pt}%      Space above
{6pt}%      Space below
{\rm}%         Body font
{}%         Indent amount (empty = no indent, \parindent = para indent)
{\bf}% Thm head font
{}%        Punctuation after thm head
{.5em}%     Space after thm head: " " = normal interword space;
{}%         Thm head spec (can be left empty, meaning `normal')

\def\Section#1{\Sec{\large #1} \setcounter{equation}{0} \vskip -6mm \indent}
\def\Sec{\@Startsection{section}{1}{\z@}
	{-3.5ex \@plus -1ex \@minus -.2ex}%
	{2.3ex \@plus.2ex}%
	{\normalfont\large\bfseries\boldmath}}
\def\@Startsection#1#2#3#4#5#6{%
	\if@noskipsec \leavevmode \fi
	\par
	\@tempskipa #4\relax
	\@afterindenttrue
	\ifdim \@tempskipa <\z@
	\@tempskipa -\@tempskipa \@afterindentfalse
	\fi
	\if@nobreak
	\everypar{}%
	\else
	\addpenalty\@secpenalty\addvspace\@tempskipa
	\fi
	\@ifstar
	{\@ssect{#3}{#4}{#5}{#6}}%
	{\@dblarg{\@Sect{#1}{#2}{#3}{#4}{#5}{#6}}}}
\def\@Sect#1#2#3#4#5#6[#7]#8{%
	\ifnum #2>\c@secnumdepth
	\let\@svsec\@empty
	\else
	\refstepcounter{#1}%
	\protected@edef\@svsec{\@seccntformat{#1}\relax}%
	\fi
	\@tempskipa #5\relax
	\ifdim \@tempskipa>\z@
	\begingroup
	#6{%
		\@hangfrom{\hskip #3\relax\@svsec \hskip -2.5mm}%
		\interlinepenalty \@M #8\@@par}
	\endgroup
	\csname #1mark\endcsname{#7}%
	\addcontentsline{toc}{#1}{%
		\ifnum #2>\c@secnumdepth \else
		\protect\numberline{\csname the#1\endcsname}%
		\fi
		#7}%
	\else
	\def\@svsechd{%
		#6{\hskip #3\relax
			\@svsec #8}%
		\csname #1mark\endcsname{#7}%
		\addcontentsline{toc}{#1}{%
			\ifnum #2>\c@secnumdepth \else
			\protect\numberline{\csname the#1\endcsname}%
			\fi
			#7}}%
	\fi
	\@xsect{#5}}
\renewenvironment{abstract}{%
	\small
	\quotation
	\noindent {\bfseries \abstractname } }%
{\if@twocolumn\else\endquotation\fi}

\def\Subsec{\@StartSubsection{subsection}{2}{\z@}%
	{-3.25ex\@plus -1ex \@minus -.2ex}%
	{1.5ex \@plus .2ex}%
	{\normalfont\normalsize\bfseries\boldmath}}
\def\@StartSubsection#1#2#3#4#5#6{%
	\if@noskipsec \leavevmode \fi
	\par
	\@tempskipa #4\relax
	\@afterindenttrue
	\ifdim \@tempskipa <\z@
	\@tempskipa -\@tempskipa \@afterindentfalse
	\fi
	\if@nobreak
	\everypar{}%
	\else
	\addpenalty\@secpenalty\addvspace\@tempskipa
	\fi
	\@ifstar
	{\@ssect{#3}{#4}{#5}{#6}}%
	{\@dblarg{\@SubSect{#1}{#2}{#3}{#4}{#5}{#6}}}}
\def\@SubSect#1#2#3#4#5#6[#7]#8{%
	\ifnum #2>\c@secnumdepth
	\let\@svsec\@empty
	\else
	\refstepcounter{#1}%
	\protected@edef\@svsec{\@seccntformat{#1}\relax}%
	\fi
	\@tempskipa #5\relax
	\ifdim \@tempskipa>\z@
	\begingroup
	#6{%
		\@hangfrom{\hskip #3\relax\@svsec\hskip -1.5mm}%
		\interlinepenalty \@M #8\@@par}
	\endgroup
	\csname #1mark\endcsname{#7}%
	\addcontentsline{toc}{#1}{%
		\ifnum #2>\c@secnumdepth \else
		\protect\numberline{\csname the#1\endcsname}%
		\fi
		#7}%
	\else
	\def\@svsechd{%
		#6{\hskip #3\relax
			\@svsec #8}%
		\csname #1mark\endcsname{#7}%
		\addcontentsline{toc}{#1}{%
			\ifnum #2>\c@secnumdepth \else
			\protect\numberline{\csname the#1\endcsname}%
			\fi
			#7}}%
	\fi
	\@xsect{#5}}
\def\list#1#2{\ifnum \@listdepth >5\relax \@toodeep \else \global
	\advance \@listdepth\@ne \fi \rightmargin \z@ \listparindent\z@
	\itemindent\z@ \csname @list\romannumeral\the\@listdepth\endcsname
	\def\@itemlabel{#1}\let\makelabel\@mklab \@nmbrlistfalse #2\relax
	\@trivlist \parskip 0pt \parindent\listparindent \advance \linewidth
	-\rightmargin \advance\linewidth -\leftmargin \advance\@totalleftmargin
	\leftmargin \parshape \@ne \@totalleftmargin \linewidth \ignorespaces}
\renewcommand{\@makecaption}[2]{\begin{center}#1. #2\end{center}}\catcode`@=12 

\theoremstyle{thmstyle}
\newtheorem{thm}{\indent Theorem}[section]
\newtheorem{lem}{\indent Lemma}[section]

\theoremstyle{remstyle}
\newtheorem{rem}{\indent \bf Remark}[section]

\newsavebox{\mygraphic}

\def\pic{\begin{picture}(0,0) \put(-210,-1250){\usebox{\mygraphic}} \end{picture}}
\newfont{\HUGEbf}{cmbx10 scaled 3500}
\definecolor{gray}{rgb}{0.9,0.9,0.9}
\def\thebibliography#1{\section*{\bf \large References}
	\list{[\arabic{enumi}]} {\settowidth \labelwidth{[#1]} \leftmargin
		\labelwidth \advance \leftmargin \labelsep \usecounter{enumi}}
	\def\newblock{\hskip .11em plus .33em minus .07em} \footnotesize \sloppy \clubpenalty
	4000 \widowpenalty 4000 \sfcode`\.=1000 \relax}

\newcommand{\eps}{\varepsilon}

\newcommand{\To}{\rightarrow}
\newcommand{\as}{{\rm d}\mathbb{P}\times{\rm d} t-a.e.}

\newcommand{\ps}{\mathbb{P}-a.s.}

\newcommand{\F}{\mathcal{F}}
\newcommand{\E}{\mathbb{E}}

\newcommand{\s}{\mathcal{S}}

\newcommand{\hcal}{\mathcal{H}}

\newcommand{\T}{[0,T]}

\newcommand{\R}{{\mathbb R}}

\newcommand{\RE}{\forall}

\newcommand {\Dis}{\displaystyle}

\def\firstpage{1}

\def\shorttitle{Multi-dimensional Mean-field Quadratic BSDEs} 
\def\shortauthor{{\it S. J. Tang\ \ and\ \ G. Yang}} 

\title{\Large \bf \boldmath\ \\ Multi-dimensional Mean-field Type Backward Stochastic Differential Equations with Diagonally Quadratic Generators$^{\ast}$}

\author{\large  Shanjian TANG$^1$\qquad Guang YANG$^{2}$}

\date{}

\begin{document}
	\maketitle
	\thispagestyle{first}
	\renewcommand{\thefootnote}{\fnsymbol{footnote}}
	\footnotetext{\hspace*{-5mm} \begin{tabular}{@{}r@{}p{13.4cm}@{}}
$^{1}$ & Department of Finance and Control Sciences, School of Mathematical Sciences, Fudan University, Shanghai 200433, China. {E-mail: sjtang@fudan.edu.cn} \\
$^2$ & Shanghai Center for Mathematical Sciences, Fudan University, Shanghai 200433, China.\\
&{E-mail: gyang19@fudan.edu.cn} \\
$^{\ast}$ & This research was supported by National Natural Science Foundation of China (Grants No. 11631004 and No. 12031009).
\end{tabular}}
\renewcommand{\thefootnote}{\arabic{footnote}}
	
\begin{abstract}
In this paper, we study the multi-dimensional backward stochastic differential equations (BSDEs) whose generator depends also on the mean of both variables. When the generator is diagonally quadratic, we prove that the BSDE admits a unique local solution with a fixed point argument. When the generator has a logarithmic growth of the off-diagonal elements (i.e., for each $i$, the $i$-th component of the generator has a logarithmic growth of the $j$-th row $z^j$ of the variable $z$ for each $j \neq i$), we give a new apriori estimate and obtain the existence and uniqueness of the global solution.
		
\vskip 4.5mm
		
		\noindent \begin{tabular}{@{}l@{ }p{10.1cm}} {\bf Keywords } &
			Multi-dimensional BSDE, Mean-field, diagonally quadratic generator, BMO martingale
		\end{tabular}
		
		\noindent {\bf 2000 MR Subject Classification } 
		60H10
	\end{abstract}
	
	\baselineskip 14pt
	
	\setlength{\parindent}{1.5em}
	
	\setcounter{section}{0}
	
\Section{Introduction} \label{section1}

In this paper, we study the existence and uniqueness of an adapted solution of the
following mean-field type BSDE:
\begin{equation}\label{eq:1.1}
Y_t=\xi+\int_t^T f(s,Y_s,\E[Y_s],Z_s,\E[Z_s]){\rm d}s-\int_t^T Z_s {\rm d}W_s, \ \ t\in\T,
\end{equation}
where $(W_t)_{t\in\T}$ is a $d$-dimensional standard Brownian motion defined on some complete probability space $(\Omega, \F, \mathbb{P})$, and $(\F_t)_{t\in\T}$ is the augmented natural filtration generated by the standard Brownian motion $W$. The terminal value $\xi$ is an $\F_T$-measurable $n$-dimensional random vector, the generator function $f(\omega, t, y, \bar y, z, \bar z):\Omega\times\T\times\R^n\times\R^n\times\R^{n\times d}\times\R^{n\times d}\To \R^n$
is $(\F_t)$-progressively measurable for each pair $(y,\bar y, z, \bar z)$, and the solution $(Y_t,Z_t)_{t\in\T}$ is a pair of $(\F_t)$-progressively measurable processes with values in $\R^n\times\R^{n\times d}$ which almost surely verifies BSDE \eqref{eq:1.1}.

When $f$ does not depend on $(\bar y, \bar z)$, BSDE~\eqref{eq:1.1} is the classical one, and it has been studied by Bismut~\cite{Bismut1973JMAA} and Pardoux and Peng~\cite{PardouxPeng1990SCL}. When the generator has a quadratic growth in the state variable $z$ and does not depend on $(\bar y, \bar z)$, BSDE~\eqref{eq:1.1} is the so-called quadratic BSDE and has been studied by Kobylanski~\cite{Kobylanski2000AP}, Briand and Hu~\cite{BriandHu2006PTRF,BriandHu2008PTRF}, Tevzadze~\cite{Tevzadze2008SPA}, Hu and Tang~\cite{HuTang2016SPA}, Xing and {\v{Z}}itkovi{\'c}~\cite{XingZitkovic2018AoP} and Fan et al.~\cite{FanHuTang2019ArXiv,FanHuTang2020ArXiv} for one-dimensional and multi-dimensional cases with bounded and unbounded terminal values.

BSDE~\eqref{eq:1.1} (so-called mean-field type BSDE) was studied by Buckdahn et al.~\cite{BuckdahnAP,BuckdahnSPA}, where they established the existence, uniqueness and a comparison theorem for the case that $f$ is uniformly Lipschitz in the last four arguments. Cheridito and Nam~\cite{cheridito} studied the existence of a class of the mean-field BSDE with quadratic growth. Carmona and Delarue~\cite{Carmona} studied some special class of quadratic forward-backward stochastic differential equations (FBSDEs) of mean-field type. Hibon, Hu, and Tang~\cite{Hibon} discussed the existence and uniqueness of one-dimensional mean-field BSDEs with quadratic growth. Hao, Wen, and Xiong~\cite{Hao} studied a class of multidimensional mean-field BSDEs with quadratic growth and small terminal value. Hao et al.~\cite{Hao2} considered the one-dimensional quadratic mean-field BSDEs when the generator depends on the laws of $(Y,Z)$. 

In this paper, we study the multidimensional mean-field BSDEs \eqref{eq:1.1} with diagonally quadratic generators and bounded terminal values by using some new methods. First, we construct a local solution with some ideas of Hu and Tang~\cite{HuTang2016SPA} and a fixed point argument. We allow the generator $f(t,y,\bar y,z,\bar z)$ to have a general growth with respect to $y$ and $\bar y$. Second, when the generator $f(t,y,\bar y,z,\bar z)$ has the additional logarithmic growth in $z$ and the additional boundedness condition with respect to $\bar z$, we build a new apriori estimate, and thus obtain the existence and uniqueness of the global solution, which is also new even when $f(t,y,\bar y,z,\bar z)$ does not depend on $(\bar y, \bar z)$.

The rest of the paper is organized as follows. In Section 2, we prepare some notations and lemmas, and state the main results of this paper. In Section 3, we prove the existence and uniqueness of the local solutions to the mean-field BSDE~\eqref{eq:1.1}. In Section 4, we give an apriori estimate and prove the existence and uniqueness of the global solutions.
	
\Section{Preliminaries and statement of main results}
\subsection{Notations}
\vspace{0.1cm}

Let $W=(W_t)_{t\geq 0}$ be a $d$-dimensional standard Brownian motion defined on a complete probability space $(\Omega, \F, \mathbb{P})$, and $(\F_t)_{t\geq 0}$ be the augmented natural filtration generated by $W$. Throughout this paper, we fix a $T \in (0, \infty)$. We endow $\Omega \times \T$ with the predictable $\sigma$-algebra $\mathcal{P}$ and $\R^n$ with its Borel $\sigma$-algebra $\mathcal{B}(\R^n)$. All the processes are assumed to be $(\F_t)_{t\in\T}$-progressively measurable, and
all equalities and inequalities between
random variables and processes are understood in the sense of $\ps$ and $\as$,  respectively. The Euclidean norm is always denoted by $|\cdot|$,  and $\|\cdot\|_{\infty}$ denotes the $L^{\infty}$-norm for one-dimensional or multidimensional random variable defined on the probability space $(\Omega, \F, \mathbb{P})$.

We define  the following four Banach spaces of stochastic processes.
By $\s^p(\R^n)$ for $p\geq 1$ , we denote the set of all $\R^n$-valued continuous adapted processes $(Y_t)_{t\in\T}$ such that
$$\|Y\|_{{\s}^p}:=\left(\E[\sup_{t\in\T} |Y_t|^p]\right)^{1/p}<+\infty.$$
By $\s^{\infty}(\R^n)$, we denote the set of all $\R^n$-valued continuous adapted processes $(Y_t)_{t\in\T}$ such that
$$\|Y\|_{{\s}^{\infty}}:=\left\|\sup_{t\in\T} |Y_t| \right\|_{\infty}<+\infty.$$
By $\hcal^p(\R^{n\times d})$ for $p\geq 1$, we denote the set of all $\R^{n\times d}$-valued $(\F_t)_{t\in\T}$-progressively measurable processes $(Z_t)_{t\in\T}$ such that
$$
\|Z\|_{\hcal^p}:=\left\{\E\left[\left(\int_0^T |Z_s|^2{\rm d}s\right)^{p/2}\right] \right\}^{1/p}<+\infty.
$$
By ${\rm BMO}(\R^{n\times d})$, we denote the set of all $Z\in \hcal^2(\R^{n\times d})$ such that
$$
\|Z\|_{\rm BMO}:=\sup_{\tau}\left\|\E_{\tau}\left[\int_{\tau}^T |Z_s|^2 {\rm d}s\right]\right\|_{\infty}^{1/2}<+\infty.
$$
Here and hereafter the supremum is taken over all $(\F_t)$-stopping times $\tau$ with values in $\T$, and $\E_{\tau}$ denotes the conditional expectation with respect to $\F_\tau$.

The spaces $\s^p_{[a,b]}(\R^n)$, $\s^{\infty}_{[a,b]}(\R^n)$, $\hcal^p_{[a,b]}(\R^{n\times d})$,  and ${\rm BMO}_{[a,b]}(\R^{n\times d})$ are identically defined for stochastic processes over the time interval $[a,b]$. We note that for $Z\in {\rm BMO}(\R^{n\times d})$, the process $\int_0^t Z_s{\rm d}W_s, t\in\T$,  is an $n$-dimensional BMO martingale.  For the theory of BMO martingales, we refer the reader to Kazamaki~\cite{Kazamaki1994book}.

For $i=1,\cdots, n$, denote by $z^i$, $y^i$, $\xi^i$ and $f^i$ the $i$th row of matrix $z\in\R^{n\times d}$,  the $i$th component of the vector $y,\xi \in \R^n$ and the generator $f$, respectively.

\subsection{Two Lemmas}

We first recall the following existence, uniqueness and an apriori estimate for one-dimensional BSDEs. The proof is given in \cite{FanHuTang2020ArXiv}.

\begin{lem}\label{lem:2.1}
We consider the following one-dimensional BSDE:
\begin{equation}\label{eq:2.1}
  Y_t=\eta+\int_t^T f(s,Z_s)\, {\rm d}s-\int_t^T Z_s\,  {\rm d}W_s, \ \ t\in\T.
\end{equation}
There exists $(U,V)\in \s^\infty(\R^n)\times {\rm BMO}(\R^{n\times d})$ such that the generator $f$  satisfies the following assumptions.

\begin{enumerate}
\renewcommand{\theenumi}{(A\arabic{enumi})}
\renewcommand{\labelenumi}{\theenumi}
\item\label{A:A1} For each $z\in \R^{1\times d}$, we have
    $$
    |f(t,z)|\leq a_t+\phi(|U_t|)+nK |V_t|^{1+\delta}+\frac{\gamma}{2} |z|^2.
    $$

\item\label{A:A2} For each $(z, \bar z)\in \R^{1\times d} \times \R^{1\times d}$, we have
    $$
     |f(t,z)-f(t,\bar z)| \leq \phi(|U_t|)\left(1+2|V_t|+|z|+|\bar z|\right) |z-\bar z|.
    $$

\item\label{A:A3} Both $|\eta|$ and $\int_0^T \alpha_t{\rm d}t$ are essentially bounded.
\end{enumerate}
Then BSDE~\eqref{eq:2.1} admits a unique solution $(Y,Z)$ such that $(Y,Z) \in \s^\infty(\R)\times {\rm BMO}(\R^{1\times d})$. Moreover, for $t\in\T$ and stopping time $\tau$ with values in $[t,T]$, we have
\begin{equation}\label{eq:2.2}
\begin{array}{lll}
|Y_t|&\leq & \Dis {1\over \gamma}\log 2+\|\eta\|_{\infty}+\Big\|\int_t^T a_s{\rm d}s\Big\|_{\infty}\vspace{0.1cm}\\
&&\Dis +\phi\left(\|U\|_{\s^{\infty}_{[t,T]}}\right)(T-t)
+\gamma^{\frac{1+\delta}{1-\delta}}C_{\delta,K,n} \|V\|_{{\rm BMO}_{[t,T]}}^{2\frac{1+\delta}{1-\delta}}(T-t),
\end{array}
\end{equation}
and
\begin{equation}\label{eq:2.3}
\begin{array}{lll}
&& \Dis\E_\tau\left[\int_\tau^T |Z_s|^2{\rm d}s\right]\vspace{0.1cm}\\
&\leq &\Dis {1\over \gamma^2} \exp(2\gamma \|\eta\|_\infty)+{1\over \gamma}\exp\left(2\gamma \Big\|\sup_{s\in [t,T]}|Y_s|\Big\|_\infty\right)\vspace{0.2cm}\\
&& \cdot \Dis\left(1+2\Big\|\int_t^T a_s {\rm d}s\Big\|_{\infty}+2\phi\left(\|U\|_{\s^{\infty}_{[t,T]}}\right)(T-t)
+2C_{\delta,K,n}\|V\|_{{\rm BMO}_{[t,T]}}^{2\frac{1+\delta}{1-\delta}}(T-t)\right),
\end{array}
\end{equation}
where
\begin{equation}\label{eq:2.4}
C_{\delta,K,n}:=\frac{1-\delta}{2}(1+\delta)^{\frac{1+\delta}{1-\delta}}
(nK)^{\frac{2}{1-\delta}}.\vspace{0.2cm}
\end{equation}
\end{lem}

Now we recall the following lemma concerning BMO-martingales.

\begin{lem}\label{lem:2.2}
For $\widetilde K>0$, there exist two constants $c_1>0$ and $c_2>0$ depending only on $\widetilde K$ such that for any BMO martingale $M$ and any one-dimensional BMO martingale $N$ such that $\|N\|_{{\rm BMO}(\mathbb P)}\le \widetilde K$, we have 
\begin{equation}\label{eq:2.5}
c_1\|M\|_{{\rm BMO}(\mathbb P)}\le \|\widetilde M\|_{{\rm BMO}(\mathbb Q)}\le c_2\|M\|_{{\rm BMO}(\mathbb P)},
\end{equation}
where ${\widetilde M}:=M-\langle M, N\rangle$ and ${\rm d} \mathbb Q:= \mathscr{E}(N)_0^\infty {\rm d} \mathbb P$.
\end{lem}

\subsection{Statement of the main results}

In this paper, we always fix several constants $\gamma > 0$, $K \geq 0$ and $\delta\in [0,1)$, a deterministic increasing continuous function $\phi(\cdot):[0,+\infty)\To [0,+\infty)$ and a $(\F_t)$-progressively measurable scalar-valued positive process $(a_t)_{t\in\T}$,.

The first main result ensures the existence and uniqueness of the local solutions for the diagonally quadratic BSDE~\eqref{eq:1.1}. We make the following assumptions.

\begin{enumerate}
\renewcommand{\theenumi}{(H\arabic{enumi})}
\renewcommand{\labelenumi}{\theenumi}
\item\label{A:H1} For $i=1,\cdots,n$ and each $(y,\bar y, z, \bar z)\in \R^n\times\R^n\times\R^{n\times d}\times\R^{n\times d}$, $f^i$ satisfies the following inequalities:
    $$
    \Dis |f^i(t,y,\bar y,z,\bar z)|\leq a_t+\phi(|y|\vee|\bar y|)+\frac{\gamma}{2} |z^i|^2+K\Big(\sum_{j\neq i} |z^j|^{1+\delta}+|\bar z|^{1+\delta}\Big).
    $$
\item\label{A:H2} For $i=1,\cdots,n$ and each $(y_k,\bar y_k, z_k, \bar z_k)\in \R^n\times\R^n\times\R^{n\times d}\times\R^{n\times d},k \in \{1,2\},$ $f^i$ satisfies the following inequalities: 
    $$
    \begin{array}{lll} 
    &&\Dis |f^i(t,y_1,\bar y_1,z_1,\bar z_1)-f^i(t,y_2,\bar y_2,z_2,\bar z_2)| \vspace{0.2cm}\\
    \Dis &\leq& \Dis \phi(|y_1|\vee|\bar y_1|\vee|y_2|\vee|\bar y_2|) \Dis \bigg[ \Big(1+|z_1|+|\bar z_1|+|z_2|+|\bar z_2|\Big)\Big(|\Delta y|+|\Delta\bar y|+|\Delta z^i|\Big) \vspace{0.2cm}\\
    \Dis&&\Dis+\Big(1+|z_1|^\delta+|\bar z_1|^{\delta}+|z_2|^\delta+|\bar z_2|^{\delta}\Big)\Big(|\Delta\bar z|+\sum_{j\neq i} |\Delta z^j|\Big)\bigg].
    \end{array}
    $$
Here~$(\Delta y,\Delta\bar y,\Delta z,\Delta\bar z)=(y_1-y_2,\bar y_1-\bar y_2, z_1-z_2, \bar z_1-\bar z_2)$.
\item\label{A:H3} There exist two positive constants $C_0$ and $C_1$ such that
$$\left\|\int_0^T a_t{\rm d}t\right\|_{\infty}\leq C_0, \ \ \|\xi\|_{\infty} \leq C_1.$$
\end{enumerate}

\begin{thm}\label{thm:2.1}
Let Assumptions \ref{A:H1}-\ref{A:H3} be satisfied. Then there exist a constant $\eps>0$ (depending only on the vector of parameters $(n,\gamma,\delta,T,K,C_0,C_1)$ and the function $\phi(\cdot)$) and a bounded subset $\mathcal{B}_\eps$ of the Banach space $\s^\infty_{[T-\eps,T]}(\R^n)\times {\rm BMO}_{[T-\eps,T]}(\R^{n\times d})$ such that BSDE~\eqref{eq:1.1} admits a unique local solution $(Y,Z)\in \mathcal{B}_\eps$ on $[T-\eps,T]$.
\end{thm}

The second main result ensures the existence and uniqueness of the global solutions for the diagonally quadratic BSDE~\eqref{eq:1.1}. The following assumptions are further required.

\begin{enumerate}
\renewcommand{\theenumi}{(H4)}
\renewcommand{\labelenumi}{\theenumi}
\item\label{A:H4} There exist a three-dimensional positive deterministic vector function $(\alpha_{t},\beta_{t},\eta_{t})_{t\in\T}$ such that for $i=1,\cdots,n$ and each $(y,\bar y, z, \bar z)\in \R^n\times\R^n\times\R^{n\times d}\times\R^{n\times d}$, $f^i$ satisfies:
    $$
    {\rm sgn}(y^i)f^i(t,y,\bar y,z,\bar z)\leq \alpha_{t}+\beta_{t}\cdot(|y|\vee |\bar y|)+\eta_{t}\log(|z|+1)+\frac{\gamma}{2} |z^i|^2.
    $$

\renewcommand{\theenumi}{(H5)}
\renewcommand{\labelenumi}{\theenumi}
\item\label{A:H5} There exists a positive constant $C_2$ such that
$$\int_0^T \Big(\alpha_t+\beta_{t}+\eta_{t}\log(1+\eta_{t})\Big){\rm d}t\leq C_2.$$

\end{enumerate}

\begin{thm}\label{thm:2.2}
Let Assumptions \ref{A:H1}-\ref{A:H5} be satisfied. Then BSDE~\eqref{eq:1.1} admits a unique global solution $(Y,Z) \in \s^\infty(\R^n)\times {\rm BMO}(\R^{n\times d})$ on $[0,T]$. 
\end{thm}

\begin{rem}\label{rem:2.1}
When the generator $f(t,y,\bar y, z, \bar z)$ does not depend on $(\bar y, \bar z)$ and BSDE~\eqref{eq:1.1} is the classical one, Theorem~\ref{thm:2.2} is still new.
\end{rem}

\Section{Local solution: the proof of Theorem~\ref{thm:2.1}}

For $ i=1, \cdots, n$, $H \in \R^{n\times d}$ and $z \in \R^{1\times d}$, denote by $H(z;i)$ the matrix in $\R^{n\times d}$ whose $i$th row is $z$ and whose $j$th row is $H^j$ for any $j\neq i$.

Let assumptions \ref{A:H1}-\ref{A:H3} be satisfied. For $(U,V) \in \s^\infty(\R^n) \times {\rm BMO}(\R^{n\times d})$,  we consider the following quadratic BSDEs:
\begin{equation}\label{eq:3.1}
Y_t^{i}=\xi^{i}+\int_t^T f^{i}(s,U_s,\E[U_s],V_s(Z_s^{i};i),\E[V_s]){\rm d}s-\int_t^T Z_s^{i} {\rm d}W_s, \ \ t\in\T, \ \ i=1,\cdots,n.
\end{equation}
From assumptions \ref{A:H1} and \ref{A:H2}, we have for each $(z, \bar z)\in \R^{1\times d} \times \R^{1\times d}$,
\begin{equation*}
\big|f^{i}(t,U_t,\E[U_t],V_t(z;i),\E[V_t])\big| \leq a_t+\phi(\|U\|_{\s^{\infty}_{[t,T]}})+\frac{\gamma}{2} |z|^2+nK|V_t|^{1+\delta}+K|\E[V_t]|^{1+\delta}; 
\end{equation*}
\begin{equation*}
\begin{array}{lll}
\Dis && \Dis \big|f^{i}(t,U_t,\E[U_t],V_t(z;i),\E[V_t])-f^{i}(t,U_t,\E[U_t],V_t(\bar z;i),\E[V_t])\big| \vspace{0.3cm}\\
\Dis &\leq& \Dis\phi(\|U\|_{\s^{\infty}_{[t,T]}})\big(1+2|V_t|+2|\E[V_t]|+|z|+|\bar z|\big)|z-\bar z|.
\end{array}
\end{equation*}
Hence the generator $f^{i}(t,U_t,\E[U_t],V_t(z;i),\E[V_t])$ satisfies the assumptions \ref{A:A1} and \ref{A:A2} in lemma~\ref{lem:2.1}. From the assumption \ref{A:H3} and lemma~\ref{lem:2.1}, we know that for each $i=1,\cdots,n$, the one-dimensional BSDE with the terminal value $\xi^i$ and the generator $f^{i}(t,U_t,\E[U_t],V_t(z;i),\E[V_t])$ admits a unique solution $(Y^i,Z^i)$ such that $(Y^i,Z^i) \in \s^\infty(\R)\times {\rm BMO}(\R^{1 \times d})$. Therefore, BSDE \eqref{eq:3.1} admits a unique solution $(Y,Z)\in \s^\infty(\R^n)\times {\rm BMO}(\R^{n\times d})$.\vspace{0.2cm}

We define the map $\Gamma$ as follows:
$$
\Gamma(U,V):=(Y,Z),\ \ \ \RE\ (U,V)\in \s^\infty(\R^n)\times {\rm BMO}(\R^{n\times d}).
$$
It is a map in the Banach space $\s^\infty(\R^n)\times {\rm BMO}(\R^{n\times d})$. From lemma~\ref{lem:2.1}, we have for $ i=1, \cdots, n$,~$t\in[0,T]$ and stopping time $\tau$ with values in $[t,T]$,
\begin{equation}\label{eq:3.2}
\begin{array}{lll}
|Y_t^i|& \leq & \Dis \frac{1}{\gamma}\log 2+\|\xi^i\|_{\infty}+\Big\|\int_t^T (a_s+K|\E[V_s]|^{1+\delta})~{\rm d}s\Big\|_{\infty}\vspace{0.2cm}\\
&&\Dis +\phi\Big(\|U\|_{\s^{\infty}_{[t,T]}}\Big)(T-t) + \gamma^{\frac{1+\delta}{1-\delta}}C_{\delta,K,n} \|V\|_{{\rm BMO}_{[t,T]}}^{2\frac{1+\delta}{1-\delta}}(T-t),
\end{array}
\end{equation}
and
\begin{equation}\label{eq:3.3}
\begin{array}{lll}
&& \Dis \E_\tau\bigg[\int_{\tau}^{T} |Z_s^i|^2~{\rm d}s\bigg]\vspace{0.2cm}\\
& \leq & \Dis \frac{1}{\gamma^2}\exp(2\gamma\|\xi^i\|_\infty)+{\frac{1}{\gamma}}\exp\Big(2\gamma \big\|\sup_{s\in [t,T]}|Y_s^i|\big\|_\infty\Big)\vspace{0.2cm}\\
&& \Dis \times \bigg[1+2\Big\|\int_t^T (a_s+K|\E[V_s]|^{1+\delta}) ~{\rm d}s\Big\|_{\infty}+2\phi\Big(\|U\|_{\s^{\infty}_{[t,T]}}\Big)(T-t) \vspace{0.2cm}\\
&&+2C_{\delta,K,n}\|V\|_{{\rm BMO}_{[t,T]}}^{2\frac{1+\delta}{1-\delta}}(T-t)\bigg].\vspace{0.1cm}
\end{array}
\end{equation}
The constant $C_{\delta,K,n}$ is defined in \eqref{eq:2.4} of lemma~\ref{lem:2.1}.
Using Young's inequality, we have
$$
\begin{array}{lll}
nK|\E[V_s]|^{1+\delta} &\leq&\Dis \frac{1}{2\|V\|_{{\rm BMO}_{[t,T]}}^2}|\E[V_s]|^2+ \frac{1-\delta}{2}(1+\delta)^{\frac{1+\delta}{1-\delta}}(nK)^{\frac{2}{1-\delta}}\|V\|_{{\rm BMO}_{[t,T]}}^{2\frac{1+\delta}{1-\delta}} \vspace{0.2cm}\\
&=& \Dis\frac{1}{2\|V\|_{{\rm BMO}_{[t,T]}}^2}|\E[V_s]|^2+ C_{\delta,K,n}\|V\|_{{\rm BMO}_{[t,T]}}^{2\frac{1+\delta}{1-\delta}}.
\end{array}
$$
Using Jensen's inequality, we have
$$
\int_t^T |\E[V_s]|^{2}{\rm d}s \leq \int_t^T \E\big[|V_s|^2\big]{\rm d}s \leq \Big\|\E_{t}\int_t^T|V_s|^2{\rm d}s\Big\|_{\infty} \leq \|V\|_{{\rm BMO}_{[t,T]}}^2.
$$
Hence we get
\begin{equation}\label{eq:3.4}
\Big\|\int_t^T nK|\E[V_s]|^{1+\delta}{\rm d}s\Big\|_{\infty} \leq \frac{1}{2}+C_{\delta,K,n}\|V\|_{{\rm BMO}_{[t,T]}}^{2\frac{1+\delta}{1-\delta}}(T-t).
\end{equation}
From \eqref{eq:3.2}, \eqref{eq:3.3}, \eqref{eq:3.4} and assumption \ref{A:H3}, we have
\begin{equation}\label{eq:3.5}
\begin{array}{lll}
\|Y\|_{\s^{\infty}_{[t,T]}} & \leq & \Dis \frac{n}{\gamma}\log 2+\frac{1}{2}+n(C_0+C_1)\vspace{0.1cm}\\
&&\Dis +n\phi\Big(\|U\|_{\s^{\infty}_{[t,T]}}\Big)(T-t)+(n\gamma^{\frac{1+\delta}{1-\delta}}+1)C_{\delta,\lambda,n} \|V\|_{{\rm BMO}_{[t,T]}}^{2\frac{1+\delta}{1-\delta}}(T-t),
\end{array}
\end{equation}
and
\begin{equation}\label{eq:3.6}
\begin{array}{lll}
\Dis \|Z\|_{{\rm BMO}_{[t,T]}}^2 & \leq &\Dis \frac{n}{\gamma^2} \exp(2\gamma C_1)+\frac{n}{\gamma}\exp\Big(2\gamma \|Y\|_{\s^{\infty}_{[t,T]}}\Big)\vspace{0.2cm}\\
&& \Dis \times \Big[2+2C_0+2\phi\Big(\|U\|_{\s^{\infty}_{[t,T]}}\Big)(T-t) + 4C_{\delta,\lambda,n}\|V\|_{{\rm BMO}_{[t,T]}}^{2\frac{1+\delta}{1-\delta}}(T-t)\Big].
\end{array}
\end{equation}
We choose
$$
K_1:=\frac{n}{\gamma}\log 2+\frac{1}{2}+n(C_0+C_1),
$$
$$
K_2:=\frac{n}{\gamma^2}\exp\big(2\gamma C_1\big)+\frac{n}{\gamma}\exp\big(4\gamma K_1\big)(2+2C_0),
$$
$$
\eps_0:=\left(\frac{K_1}{n\phi\left(2K_1\right)
+(n\gamma^{\frac{1+\delta}{1-\delta}}+1)C_{\delta,\lambda,n} (2K_2)^{\frac{1+\delta}{1-\delta}}}\right) \bigwedge \left(\frac{{\gamma\over n}\exp\left(-4\gamma K_1\right)K_2}{2\phi\left(2K_1\right)
+4C_{\delta,\lambda,n} (2K_2)^{\frac{1+\delta}{1-\delta}}}\right)>0.\vspace{0.3cm}
$$
From \eqref{eq:3.5} and \eqref{eq:3.6}, we have for $\eps \in (0,\eps_0]$, if
$$
\|U\|_{\s^{\infty}_{[T-\eps,T]}} \leq 2K_1, \ \ \ \|V\|_{{\rm BMO}_{[T-\eps,T]}}^2\leq 2K_2,
$$
then
$$
\|Y\|_{\s^{\infty}_{[T-\eps,T]}} \leq 2K_1, \ \ \ \|Z\|_{{\rm BMO}_{[T-\eps,T]}}^2 \leq 2K_2.
$$
Define
\begin{equation}\label{eq:3.7}
\begin{array}{l}
\Dis\mathcal{B}_\eps:=\left\{(U,V)\in \s^\infty(\R^n)\times {\rm BMO}(\R^{n\times d}): \|U\|_{\s^{\infty}_{[T-\eps,T]}}\leq 2K_1,\ \|V\|_{{\rm BMO}_{[T-\eps,T]}}^2\leq 2K_2\right\}.\vspace{0.2cm}
\end{array}
\end{equation}
Then for each $\eps \in (0,\eps_0]$, $\Gamma$ maps the Banach space $\mathcal{B}_\eps$ to itself.

It remains to prove that there exists a real constant $\eps \in (0,\eps_0]$ depending only on constants $(n,\gamma,\delta,T,K,C_0,C_1)$ and the function $\phi(\cdot)$ such that $\Gamma$ is a contraction in $\mathcal{B}_\eps$. For a given $\eps \in (0,\eps_0]$ and $(U,V),(\widetilde U,\widetilde V) \in \mathcal{B}_\eps$, we denote
$$
(Y,Z):=\Gamma(U,V), \ \ \ (\widetilde Y,\widetilde Z):=\Gamma(\widetilde U,\widetilde V).
$$
Then we have for $i=1, \cdots, n$ and $t \in [T-\eps,T]$,
\begin{equation*}
\begin{array}{ll}
\Dis Y_t^{i}=\xi^{i}+\int_t^T f^{i}(s,U_s,\E[U_s],V_s(Z_s^{i};i),\E[V_s]){\rm d}s-\int_t^T Z_s^{i} {\rm d}W_s,\vspace{0.2cm}\\
\Dis \widetilde Y_t^{i}=\xi^{i}+\int_t^T f^{i}(s,\widetilde U_s,\E[\widetilde U_s],\widetilde V_s(\widetilde Z_s^{i};i),\E[\widetilde V_s]){\rm d}s-\int_t^T \widetilde Z_s^{i} {\rm d}W_s.
\end{array}
\end{equation*}
For $i=1, \cdots, n$ and $s \in [T-\eps,T]$, define
\begin{equation*}
\begin{array}{ll}
\Dis \Delta_s^{1,i}:=f^{i}(s,U_s,\E[U_s],V_s(Z_s^{i};i),\E[V_s])-f^{i}(s,U_s,\E[U_s],V_s(\widetilde Z_s^{i};i),\E[V_s]),\vspace{0.2cm}\\
\Dis \Delta_s^{2,i}:=f^{i}(s,U_s,\E[U_s],V_s(\widetilde Z_s^{i};i),\E[V_s])-f^{i}(s,\widetilde U_s,\E[\widetilde U_s],V_s(\widetilde Z_s^{i};i),\E[V_s]),\vspace{0.2cm}\\
\Dis \Delta_s^{3,i}:=f^{i}(s,\widetilde U_s,\E[\widetilde U_s],V_s(\widetilde Z_s^{i};i),\E[V_s])-f^{i}(s,\widetilde U_s,\E[\widetilde U_s],\widetilde V_s(\widetilde Z_s^{i};i),\E[\widetilde V_s]),\vspace{0.2cm}\\
\Dis (\Delta Y, \Delta Z, \Delta U, \Delta V, \Delta \E[V]):=(Y-\widetilde Y, Z-\widetilde Z, U-\widetilde U, V-\widetilde V, \E[V]-\E[\widetilde V]).
\end{array}
\end{equation*}
Then we have
\begin{equation}\label{eq:3.8}
\Delta Y_t^{i}+\int_t^T \Delta Z_s^{i} {\rm d}W_s-\int_t^T \Delta_s^{1,i}{\rm d}s=\int_t^T (\Delta_s^{2,i}+\Delta_s^{3,i}){\rm d}s, \ \ t \in [T-\eps,T].
\end{equation}
From assumption \ref{A:H2} and $(U,V),(\widetilde U,\widetilde V) \in \mathcal{B}_\eps$, for $i=1, \cdots, n$ and $s \in [T-\eps,T]$, we have
\begin{equation}\label{eq:3.9}
|\Delta_s^{1,i}| \leq \phi(2K_1)\big(1+2|V_s|+2|\E[V_s]|+|Z_s|+|\widetilde Z_s|\big)|\Delta Z_s^{i}|,
\end{equation}
and
\begin{equation}\label{eq:3.10}
\begin{array}{ll}
\Dis |\Delta_s^{2,i}| \leq 2\phi(2K_1)\big(1+2|V_s|+2|\E[V_s]|+2|\widetilde Z_s|\big)\|\Delta U\|_{\s^{\infty}};\vspace{0.2cm}\\
|\Delta_s^{3,i}| \leq \phi(2K_1)\big(1+|V_s|^{\delta}+|\E[V_s]|^{\delta}+|\widetilde V_s|^{\delta}+|\E[\widetilde V_s]|^{\delta}+2|\widetilde Z_s|^{\delta}\big)(|\Delta \E[V_s]|+\sqrt{n}|\Delta V_s|).
\end{array}
\end{equation}
For $i=1, \cdots, n$ and $s\in [T-\eps,T]$, from \eqref{eq:3.9} we know there exists a $\R^{d}$-valued process $\Lambda(i)$ such that
\begin{equation}\label{eq:3.11}
\Delta_s^{1,i} = \Delta Z_s^{i} \cdot \Lambda_{s}(i), \ \ \ |\Lambda_{s}(i)| \leq \phi(2K_1)\big(1+2|V_s|+2|\E[V_s]|+|Z_s|+|\widetilde Z_s|\big).
\end{equation}
We can take $\Lambda_s(i)=0$ on $[0,T-\eps]$. Then $\widetilde W_t(i):=W_t-\int_0^t \Lambda_s(i) {\rm d}s $ is a Brownian motion under the equivalent probability measure $\mathbb Q^{i}$ defined by
$$
{\rm d} \mathbb Q^{i}:= \mathscr{E}(\Lambda(i) \cdot W)_0^T {\rm d} \mathbb P.
$$
From \eqref{eq:3.11}, we have
\begin{equation*}
\begin{array}{lll}
\Dis\E_{t}\int_t^T |\Lambda_s(i)|^2 {\rm d}s &\leq& \Dis [\phi(2K_1)]^{2}\E_{t}\int_t^T \big(1+2|V_s|+2|\E[V_s]|+|Z_s|+|\widetilde Z_s|\big)^2 {\rm d}s \vspace{0.2cm}\\
&\leq& \Dis 5[\phi(2K_1)]^{2}\E_{t}\int_t^T \big(1+4|V_s|^2+4|\E[V_s]|^2+|Z_s|^2+|\widetilde Z_s|^2\big) {\rm d}s\vspace{0.2cm}\\
&\leq& \Dis 5[\phi(2K_1)]^{2}(T+12K_2+\E\int_t^T 4|V_s|^2 {\rm d}s)\vspace{0.2cm}\\
&\leq& 5(T+20K_2)[\phi(2K_1)]^{2}.
\end{array}
\end{equation*}
Therefore $\|\Lambda(i) \cdot W\|_{{\rm BMO}(\mathbb P)}^{2}\le \widetilde K^2:=5(T+20K_2)\big[\phi(2K_1)\big]^{2}$. From \eqref{eq:3.8} and \eqref{eq:3.11}, we have
\begin{equation}\label{eq:3.12}
\Delta Y_t^{i}+\int_t^T \Delta Z_s^{i} {\rm d}\widetilde W_s(i)=\int_t^T (\Delta_s^{2,i}+\Delta_s^{3,i}){\rm d}s, \ \ t \in [T-\eps,T].
\end{equation}
Taking square and the conditional expectation with respect to $\mathbb Q^{i}$, we have for $t \in [T-\eps,T]$, 
\begin{equation}\label{eq:3.13}
\begin{array}{lll}
\Dis |\Delta Y_t^{i}|^2+\E_{t}^{\mathbb Q^{i}}\Big[\int_t^T |\Delta Z_s^{i}|^2 {\rm d}s\Big] & = & \Dis \E_{t}^{\mathbb Q^{i}}\Big[\big(\int_t^T (\Delta_s^{2,i}+\Delta_s^{3,i}){\rm d}s\big)^2\Big]\vspace{0.2cm}\\
& \leq & \Dis 2\E_{t}^{\mathbb Q^{i}}\Big[\Big(\int_t^T |\Delta_s^{2,i}| {\rm d}s\Big)^{2}\Big]+2\E_{t}^{\mathbb Q^{i}}\Big[\Big(\int_t^T |\Delta_s^{3,i}| {\rm d}s\Big)^2\Big].
\end{array}
\end{equation}
From \eqref{eq:3.10} and Cauchy-Schwarz inequality, we obtain that for $t \in [T-\eps,T]$, 
\begin{align}
&~\Dis\E_{t}^{\mathbb Q^{i}}\Big[\Big(\int_t^T |\Delta_s^{2,i}| {\rm d}s\Big)^{2}\Big]\vspace{0.2cm}\notag\\
\leq&\Dis~4\big[\phi(2K_1)\big]^{2}\big\|\Delta U\big\|_{\s^{\infty}}^{2} \E_{t}^{\mathbb Q^{i}}\Big[\Big(\int_t^T \big(1+2|V_s|+2|\E[V_s]|+2|\widetilde Z_s|\big) {\rm d}s\Big)^{2}\Big]\vspace{0.2cm}\notag\\
\leq&~\Dis4(T-t)\big[\phi(2K_1)\big]^{2}\big\|\Delta U\big\|_{\s^{\infty}}^{2} \E_{t}^{\mathbb Q^{i}}\Big[\int_t^T \big(1+2|V_s|+2|\E[V_s]|+2|\widetilde Z_s|\big)^2 {\rm d}s\Big]\vspace{0.2cm}\notag\\
\leq&~\Dis16(T-t)\big[\phi(2K_1)\big]^{2}\big\|\Delta U\big\|_{\s^{\infty}}^{2} \E_{t}^{\mathbb Q^{i}}\Big[\int_t^T \big(1+4|V_s|^{2}+4|\E[V_s]|^{2}+4|\widetilde Z_s|^{2}\big) {\rm d}s\Big]\vspace{0.2cm}\notag\\
\leq&~\Dis16\eps\big[\phi(2K_1)\big]^{2}\big\|\Delta U\big\|_{\s^{\infty}}^{2}\Big(T+4\|V\|_{{\rm BMO}(\mathbb Q^{i})}^{2}+4\|\widetilde Z\|_{{\rm BMO}(\mathbb Q^{i})}^{2}+4\E\int_{t}^{T}|V_s|^{2} {\rm d}s\Big)\vspace{0.2cm}\notag\\
\leq&~\Dis16\eps\big[\phi(2K_1)\big]^{2}\big\|\Delta U\big\|_{\s^{\infty}}^{2}\Big(T+4\|V\|_{{\rm BMO}(\mathbb Q^{i})}^{2}+4\|\widetilde Z\|_{{\rm BMO}(\mathbb Q^{i})}^{2}+4\|V\|_{{\rm BMO}(\mathbb P)}^{2}\Big),
\end{align}
and
\begin{equation}\label{eq:3.15}
\begin{array}{lll}
&&\Dis \E_{t}^{\mathbb Q^{i}}\Big[\big(\int_t^T |\Delta_s^{3,i}| {\rm d}s\big)^{2}\Big]\vspace{0.2cm}\\
&\leq&\Dis\big[\phi(2K_1)\big]^{2}\E_{t}^{\mathbb Q^{i}}\Big[\Big(\int_t^T \big(1+|V_s|^{\delta}+|\E[V_s]|^{\delta}+|\widetilde V_s|^{\delta}+|\E[\widetilde V_s]|^{\delta}+2|\widetilde Z_s|^{\delta}\big)\vspace{0.2cm}\\
&&\Dis\times\big(|\Delta \E[V_s]|+\sqrt{n}|\Delta V_s|\big) {\rm d}s\Big)^{2}\Big]\vspace{0.2cm}\\
&\leq&\Dis\big[\phi(2K_1)\big]^{2}\E_{t}^{\mathbb Q^{i}}\Big[\int_t^T \big(1+|V_s|^{\delta}+|\E[V_s]|^{\delta}+|\widetilde V_s|^{\delta}+|\E[\widetilde V_s]|^{\delta}+2|\widetilde Z_s|^{\delta}\big)^{2} {\rm d}s\vspace{0.2cm}\\
&& \times\Dis\int_t^T \big(|\Delta \E[V_s]|+\sqrt{n}|\Delta V_s|\big)^{2} {\rm d}s\Big]\vspace{0.2cm}\\
&\leq&\Dis\big[\phi(2K_1)\big]^{2}\E_{t}^{\mathbb Q^{i}}\bigg[\Big(\int_t^T \big(1+|V_s|^{\delta}+|\E[V_s]|^{\delta}+|\widetilde V_s|^{\delta}+|\E[\widetilde V_s]|^{\delta}+2|\widetilde Z_s|^{\delta}\big)^{2} {\rm d}s\Big)^{2}\bigg]^{\frac{1}{2}}\vspace{0.2cm}\\
&&\times\Dis\E_{t}^{\mathbb Q^{i}}\bigg[\Big(\int_t^T (|\Delta \E[V_s]|+\sqrt{n}|\Delta V_s|)^{2} {\rm d}s\Big)^2\bigg]^{\frac{1}{2}}.
\end{array}
\end{equation}
Let $L_4$ be a constant such that for any BMO martingale $M$, we have
$$
\sup_{\tau}\E_{\tau}[(\langle M \rangle)^2)] \ \leq L_{4}^{4} \|M\|_{\rm BMO}^{4}.
$$
Using Cauchy-Schwarz inequality and H\"{o}lder's inequality, we get
\begin{align}
&\Dis~\E_{t}^{\mathbb Q^{i}}\bigg[\Big(\int_t^T \big(1+|V_s|^{\delta}+|\E[V_s]|^{\delta}+|\widetilde V_s|^{\delta}+|\E[\widetilde V_s]|^{\delta}+2|\widetilde Z_s|^{\delta}\big)^{2} {\rm d}s\Big)^{2}\bigg]^{\frac{1}{2}}\vspace{0.2cm}\notag\\
\leq&\Dis~\E_{t}^{\mathbb Q^{i}}\bigg[\Big(\int_t^T 6\big(1+|V_s|^{2\delta}+|\E[V_s]|^{2\delta}+|\widetilde V_s|^{2\delta}+|\E[\widetilde V_s]|^{2\delta}+4|\widetilde Z_s|^{2\delta}\big) {\rm d}s\Big)^{2}\bigg]^{\frac{1}{2}}\vspace{0.2cm}\notag\\
\leq&\Dis ~6\E_{t}^{\mathbb Q^{i}}\bigg[\Big[\eps+\eps^{1-\delta}\Big(\int_t^T|V_s|^{2}{\rm d}s\Big)^{\delta}+\eps^{1-\delta}\Big(\int_t^T|\E[V_s]|^{2}{\rm d}s\Big)^{\delta}+\eps^{1-\delta}\Big(\int_t^T|\widetilde V_s|^{2}{\rm d}s\Big)^{\delta}\vspace{0.2cm}\notag\\
&\Dis ~+\eps^{1-\delta}\Big(\int_t^T|\E[\widetilde V_s]|^{2}{\rm d}s\Big)^{\delta}+4\eps^{1-\delta}\Big(\int_t^T|\widetilde Z_s|^{2}{\rm d}s\Big)^{\delta}\Big]^{2}\bigg]^{\frac{1}{2}} \vspace{0.2cm}\notag\\
=&\Dis ~6\eps^{1-\delta} \E_{t}^{\mathbb Q^{i}}\bigg[\Big[\eps^{\delta}+\Big(\int_t^T|V_s|^{2}{\rm d}s\Big)^{\delta}+\Big(\int_t^T|\E[V_s]|^{2}{\rm d}s\Big)^{\delta}+\Big(\int_t^T|\widetilde V_s|^{2}{\rm d}s\Big)^{\delta}\vspace{0.2cm}\notag\\
&\Dis~+\Big(\int_t^T|\E[\widetilde V_s]|^{2}{\rm d}s\Big)^{\delta}+4\Big(\int_t^T|\widetilde Z_s|^{2}{\rm d}s\Big)^{\delta}\Big]^{2}\bigg]^{\frac{1}{2}}\vspace{0.2cm}\notag\\
\leq&\Dis ~6\eps^{1-\delta} \E_{t}^{\mathbb Q^{i}}\bigg[\Big(T^{\delta}+\int_t^T|V_s|^{2}{\rm d}s+1+\int_t^T|\E[V_s]|^{2}{\rm d}s+1+\int_t^T|\widetilde V_s|^{2}{\rm d}s+1\vspace{0.2cm}\notag\\
&\Dis~+\int_t^T|\E[\widetilde V_s]|^{2}{\rm d}s+1+4\int_t^T|\widetilde Z_s|^{2}{\rm d}s+4\Big)^{2}\bigg]^{\frac{1}{2}} \vspace{0.2cm}\notag\\
\leq&\Dis ~6\sqrt{6}\eps^{1-\delta} \E_{t}^{\mathbb Q^{i}}\bigg[\Big(T^{\delta}+8\Big)^2+\Big(\int_t^T|V_s|^{2}{\rm d}s\Big)^{2}+\Big(\int_t^T|\E[V_s]|^{2}{\rm d}s\Big)^{2}+\Big(\int_t^T|\widetilde V_s|^{2}{\rm d}s\Big)^{2}\vspace{0.2cm}\notag\\
&\Dis~+\Big(\int_t^T|\E[\widetilde V_s]|^{2}{\rm d}s\Big)^{2}+16\Big(\int_t^T|\widetilde Z_s|^{2}{\rm d}s\Big)^{2}\bigg]^{\frac{1}{2}} \vspace{0.2cm}\notag\\
\leq&\Dis ~6\sqrt{6}\eps^{1-\delta}\bigg[T^{\delta}+8+\E_{t}^{\mathbb Q^{i}}\Big[\Big(\int_t^T|V_s|^{2}{\rm d}s\Big)^{2}\Big]^{\frac{1}{2}}+\int_t^T|\E[V_s]|^{2}{\rm d}s+\E_{t}^{\mathbb Q^{i}}\Big[\Big(\int_t^T|\widetilde V_s|^{2}{\rm d}s\Big)^{2}\Big]^{\frac{1}{2}}\vspace{0.2cm}\notag\\
&\Dis~+\int_t^T|\E[\widetilde V_s]|^{2}{\rm d}s+4\E_{t}^{\mathbb Q^{i}}\Big[\Big(\int_t^T|\widetilde Z_s|^{2}{\rm d}s\Big)^{2}\Big]^{\frac{1}{2}}\bigg] \vspace{0.2cm}\notag\\
\leq&\Dis ~6\sqrt{6}\eps^{1-\delta}\bigg[L_4^{2}\Big(\|V\|_{{\rm BMO}(\mathbb Q^{i})}^{2}+\|\widetilde V\|_{{\rm BMO}(\mathbb Q^{i})}^{2}\Big)+4L_4^{2}\big\|\widetilde Z\big\|_{{\rm BMO}(\mathbb Q^{i})}^{2}+\int_t^T \E\big[|V_s|^{2}+|\widetilde V_s|^{2}\big]{\rm d}s\bigg] \vspace{0.2cm}\notag\\
&\Dis~+6\sqrt{6}\eps^{1-\delta}\Big(T^{\delta}+8\Big)\vspace{0.2cm}\notag\\
\leq&\Dis ~6\sqrt{6}\eps^{1-\delta}\bigg[L_4^{2}\Big(\|V\|_{{\rm BMO}(\mathbb Q^{i})}^{2}+\|\widetilde V\|_{{\rm BMO}(\mathbb Q^{i})}^{2}\Big)+4L_4^{2}\big\|\widetilde Z\big\|_{{\rm BMO}(\mathbb Q^{i})}^{2}+\big\|V\big\|_{{\rm BMO}(\mathbb P)}^{2}+\big\|\widetilde V\big\|_{{\rm BMO}(\mathbb P)}^{2}\bigg]\vspace{0.2cm}\notag\\
&\Dis~+6\sqrt{6}\eps^{1-\delta}\Big(T^{\delta}+8\Big),
\end{align}
and
\begin{equation}\label{eq:3.17}
\begin{array}{lll}
&&\Dis\E_{t}^{\mathbb Q^{i}}\bigg[\Big(\int_t^T \big(|\Delta \E[V_s]|+\sqrt{n}|\Delta V_s|\big)^{2} {\rm d}s\Big)^2\bigg]^{\frac{1}{2}}\vspace{0.2cm}\\
&\leq&\Dis\E_{t}^{\mathbb Q^{i}}\bigg[\Big(\int_t^T 2\big(|\Delta \E[V_s]|^{2}+n|\Delta V_s|^{2}\big) {\rm d}s\Big)^2\bigg]^{\frac{1}{2}}\vspace{0.2cm}\\
&\leq&\Dis2 \E_{t}^{\mathbb Q^{i}}\bigg[2\Big(\int_t^T |\Delta \E[V_s]|^{2} {\rm d}s\Big)^{2} + 2\Big(\int_t^T n|\Delta V_s|^{2} {\rm d}s\Big)^{2}\bigg]^{\frac{1}{2}}\vspace{0.2cm}\\
&\leq&\Dis2\sqrt{2}\int_t^T |\Delta \E[V_s]|^{2} {\rm d}s + 2\sqrt{2}n\E_{t}^{\mathbb Q^{i}}\bigg[\Big(\int_t^T |\Delta V_s|^{2} {\rm d}s\Big)^{2}\bigg]^{\frac{1}{2}}\vspace{0.2cm}\\
&\leq&\Dis2\sqrt{2}\int_t^T \E\big[|\Delta V_s|^{2}\big] {\rm d}s + 2\sqrt{2}nL_4^{2}\big\|\Delta V\big\|_{{\rm BMO}(\mathbb Q^{i})}^{2}\vspace{0.2cm}\\
&\leq&\Dis2\sqrt{2}\big\|\Delta V\big\|_{{\rm BMO}(\mathbb P)}^{2} + 2\sqrt{2}nL_4^{2}\big\|\Delta V\big\|_{{\rm BMO}(\mathbb Q^{i})}^{2}.\vspace{0.1cm}
\end{array}
\end{equation}
Comining \eqref{eq:3.13}-\eqref{eq:3.17} and lemma~\ref{lem:2.2}, we know there exist two positive constants $c_1$ and $c_2$ depending only on $\widetilde K^2=5(T+20K_2)\big[\phi(2K_1)\big]^2$ such that for $t \in [T-\eps,T]$,
\begin{equation}\label{eq:3.18}
\begin{array}{lll} 
&&\Dis |\Delta Y_t^{i}|^2+\E_{t}^{\mathbb Q^{i}}\Big[\int_t^T |\Delta Z_s^{i}|^2 {\rm d}s\Big] \vspace{0.2cm}\\
&\leq&\Dis 32\eps\big[\phi(2K_1)\big]^{2}\Big(T+16c_2^{2}K_2+8K_2\Big) \big\|\Delta U\big\|_{\s^{\infty}}^{2}\vspace{0.2cm}\\
&&\Dis+48\sqrt{3}\eps^{1-\delta}\big[\phi(2K_1)\big]^{2}\Big(1+nL_4^{2}c_2^{2}\Big)\Big(T^{\delta}+8+12L_4^{2}c_2^{2}K_2+4K_2\Big)\big\|\Delta V\big\|_{{\rm BMO}(\mathbb P)}^{2}.
\end{array}
\end{equation}
Therefore, we have on the interval $[T-\eps, T]$,
\begin{align}\label{eq:3.19} 
&\Dis ~\big\|\Delta Y\big\|_{\s^{\infty}}^2+c_1^{2}\big\|\Delta Z\big\|_{{\rm BMO}(\mathbb P)}^2 \vspace{0.2cm}\notag\\
\leq&\Dis ~64n\eps\big[\phi(2K_1)\big]^{2}\Big(T+16c_2^{2}K_2+8K_2\Big) \big\|\Delta U\big\|_{\s^{\infty}}^{2}\vspace{0.2cm}\notag\\
&\Dis~+96\sqrt{3}n\eps^{1-\delta}\big[\phi(2K_1)\big]^{2}\Big(1+nL_4^{2}c_2^{2}\Big)\Big(T^{\delta}+8+12L_4^{2}c_2^{2}K_2+4K_2\Big)\big\|\Delta V\big\|_{{\rm BMO}(\mathbb P)}^{2}.
\end{align}

Hence $\Gamma$ is a contraction in $\mathcal{B}_\eps$ for sufficiently small $\eps$. From fixed-point Theorem, we get Theorem~\ref{thm:2.1}. The proof is complete.

\Section{Global solution: the proof of Theorem~\ref{thm:2.2}}

We first give an inequality.

\begin{lem}\label{lem:4.1}
For $x,y,C >0$, we have
\begin{equation}\label{eq:4.1}
C\log(1+x) \leq \frac{x^{2}}{y}+C\log(1+Cy). \vspace{0.2cm}
\end{equation}
\end{lem}

{\bf Proof}\ \ When $0<x\leq Cy$, we have
$$
C\log(1+x) \leq C\log(1+Cy) \leq \frac{x^{2}}{y}+C\log(1+Cy).
$$
When $x > Cy$, we have 
\begin{equation*}
C\log(1+x) \leq Cx \leq \frac{x^{2}}{y} \leq \frac{x^{2}}{y}+C\log(1+Cy).
\end{equation*}
The proof is complete.\vspace{0.2cm}

Now we give an apriori estimate.
\begin{lem}\label{lem:4.2}
Let Assumptions \ref{A:H1}-\ref{A:H5} be satisfied, $(Y,Z) \in \s^{\infty}_{[t_0,T]}(\R^n) \times {\hcal^2}_{[t_0,T]}(\R^{n\times d})$ is a solution of BSDE \eqref{eq:1.1} on $[t_0,T]$, then there exist a positive constant $\lambda$, depending only on the vector of parameters $(n,\gamma,T,C_1,C_2)$, such that

\begin{equation}\label{eq:4.2}
\|Y\|_{\s^{\infty}_{[t_0,T]}} \leq \lambda.\vspace{0.2cm}
\end{equation}
\end{lem} 

{\bf Proof}\ \ Define
$$
u(x)=\gamma^{-2}\big(\exp(\gamma|x|)-\gamma|x|-1\big), \ \ x \in \R.
$$
Then we have for $x \in \R$,
$$
u'(x)=\gamma^{-1}\big(\exp(\gamma|x|)-1\big){\rm sgn}(x), \ \ u''(x)=\exp(\gamma|x|), \ \ u''(x)-\gamma |u'(x)| =1.
$$
Using It\^{o}'s formula to compute $u(Y_{t}^{i})$, we have for $t\in [t_0,T]$,
\begin{align}\label{eq:4.3}
\Dis u(Y_{t}^{i}) = &~\Dis u(\xi^{i})-\int_{t}^{T}u'(Y_{s}^{i})Z_{s}^{i}{\rm d}W_s\vspace{0.2cm}\notag\\
&\Dis ~+\int_{t}^{T}\Big[u'(Y_{s}^{i})f^{i}(s,Y_{s},\E(Y_{s}),Z_{s},\E(Z_{s}))-\frac{1}{2}u''(Y_{s}^{i})|Z_{s}^{i}|^{2}\Big] {\rm d}s \vspace{0.2cm}\notag\\
\leq & \Dis ~u(\xi^{i})-\int_{t}^{T}u'(Y_{s}^{i})Z_{s}^{i}{\rm d}W_s-\frac{1}{2}\int_{t}^{T}\exp(\gamma|Y_{s}^{i}|)|Z_{s}^{i}|^{2} {\rm d}s \vspace{0.2cm}\notag\\
& \Dis ~+\int_{t}^{T}\gamma^{-1}\Big[\exp(\gamma|Y_{s}^{i}|)-1\Big]\Big[\frac{\gamma}{2}|Z_{s}^{i}|^{2}+\alpha_{s}+\beta_{s}(|Y_{s}|\vee|\E(Y_{s})|)+\eta_{s}\log(|Z_{s}|+1)\Big] {\rm d}s \vspace{0.2cm}\notag\\
= & \Dis ~u(\xi^{i})-\int_{t}^{T}u'(Y_{s}^{i})Z_{s}^{i}{\rm d}W_s-\frac{1}{2}\int_{t}^{T}|Z_{s}^{i}|^{2} {\rm d}s \vspace{0.2cm}\notag\\
& \Dis ~+\int_{t}^{T}\gamma^{-1}\Big[\exp(\gamma|Y_{s}^{i}|)-1\Big]\Big[\alpha_{s}+\beta_{s}(|Y_{s}|\vee|\E(Y_{s})|)+\eta_{s}\log(|Z_{s}|+1)\Big] {\rm d}s.
\end{align}

Using (\ref{eq:4.1}) and taking
\begin{equation*}
x=|Z_{s}|, \ \ y=\frac{2n}{\gamma}\exp(\gamma\|Y\|_{\s^{\infty}_{[s,T]}}), \ \ C=\eta_{s},
\end{equation*}
we have
\begin{equation}\label{eq:4.4}
\begin{array}{lll}
\Dis && \Dis \eta_{s}\log(|Z_{s}|+1) \vspace{0.2cm}\\
& \leq & \Dis \frac{\gamma}{2n}\exp(-\gamma\|Y\|_{\s^{\infty}_{[s,T]}})|Z_{s}|^2+\eta_{s}\log\big(1+\frac{2}{\gamma}n \eta_{s}\exp(\gamma\|Y\|_{\s^{\infty}_{[s,T]}})\big) \vspace{0.2cm}\\
& \leq & \Dis \frac{\gamma}{2n}\exp(-\gamma\|Y\|_{\s^{\infty}_{[s,T]}})|Z_{s}|^2+\eta_{s}\log\big(1+(\frac{2}{\gamma}n \eta_{s}+1)\exp(\gamma\|Y\|_{\s^{\infty}_{[s,T]}})\big).
\end{array}
\end{equation}
It is easy to check that
\begin{equation*}
\begin{array}{lll}
\Dis \log(1+x) \leq 1+\log x, \ \ \forall x \geq 1; \vspace{0.2cm}\\
\Dis \log(1+kx) \leq k+\log(1+x), \ \ \forall k \geq 0, \ \ x \geq 0. 
\end{array}
\end{equation*}
Hence we have
\begin{equation*}
\begin{array}{lll}
& &\Dis \log\big(1+(\frac{2}{\gamma}n \eta_{s}+1)\exp(\gamma\|Y\|_{\s^{\infty}_{[s,T]}})\big) \\
&\leq & \Dis 1+\gamma\|Y\|_{\s^{\infty}_{[s,T]}}+\log(\frac{2}{\gamma}n \eta_{s}+1) \\
& \leq & \Dis 1+\gamma\|Y\|_{\s^{\infty}_{[s,T]}}+\log(\eta_{s}+1)+\frac{2n}{\gamma}.
\end{array}
\end{equation*}
From (\ref{eq:4.4}), we have 
\begin{equation}\label{eq:4.5}
\begin{array}{lll}
\Dis && \Dis \eta_{s}\log(|Z_{s}|+1) \vspace{0.2cm}\\
& \leq & \Dis \frac{\gamma}{2n}\exp(-\gamma\|Y\|_{\s^{\infty}_{[s,T]}})|Z_{s}|^2+\eta_{s}\Big(1+\gamma\|Y\|_{\s^{\infty}_{[s,T]}}+\log(\eta_{s}+1)+\frac{2n}{\gamma}\Big).
\end{array}
\end{equation}
Let 
\begin{equation*}
k_s=\eta_{s}\big(1+\log(\eta_{s}+1)+\frac{2n}{\gamma}\big).
\end{equation*}
From (\ref{eq:4.3}) and (\ref{eq:4.5}), we have
\begin{equation}\label{eq:4.6}
\begin{array}{lll}
\Dis u(Y_{t}^{i}) & \leq & \Dis  u(\xi^{i})-\int_{t}^{T}u'(Y_{s}^{i})Z_{s}^{i}{\rm d}W_s+\frac{1}{2}\int_{t}^{T}\big(-|Z_{s}^{i}|^{2}+\frac{1}{n}|Z_{s}|^{2}\big) {\rm d}s \vspace{0.2cm}\\
& & \Dis +\int_{t}^{T}\gamma^{-1}\exp(\gamma|Y_{s}^{i}|)\Big(\alpha_{s}+(\beta_{s}+\gamma\eta_{s})\|Y\|_{\s^{\infty}_{[s,T]}}+k_s\Big) {\rm d}s, \ \ t\in [t_0,T]
\end{array}
\end{equation}
Hence, we have
\begin{equation}\label{eq:4.7}
\begin{array}{lll}
\Dis \sum_{i=1}^{n}u(Y_{t}^{i}) & \leq & \Dis  \sum_{i=1}^{n}u(\xi^{i})-\int_{t}^{T} \sum_{i=1}^{n}u'(Y_{s}^{i})Z_{s}^{i}{\rm d}W_s \vspace{0.2cm}\\
& & \Dis +\int_{t}^{T}\gamma^{-1}\Big(\alpha_{s}+(\beta_{s}+\gamma\eta_{s})\|Y\|_{\s^{\infty}_{[s,T]}}+k_s\Big)\sum_{i=1}^{n}\exp(\gamma|Y_{s}^{i}|) {\rm d}s, \ \ t\in [t_0,T].
\end{array}
\end{equation} 
Noting that 
$$
\frac{1}{2}\gamma^{-2}\big(\exp(\gamma|x|)-2\big) \leq u(x) \leq \gamma^{-2}\exp(\gamma|x|),
$$ 
we have for $t\in [t_0,T]$,
\begin{equation}\label{eq:4.8}
\begin{array}{lll}
\Dis \frac{1}{2}\gamma^{-2}\sum_{i=1}^{n}\big[\exp(\gamma|Y_{t}^{i}|)-2\big] & \leq & \Dis  \gamma^{-2}n\exp(\gamma\|\xi\|_{\infty})-\int_{t}^{T} \sum_{i=1}^{n}u'(Y_{s}^{i})Z_{s}^{i}{\rm d}W_s \vspace{0.2cm}\\
& & \Dis +\gamma^{-1}\int_{t}^{T}\Big(\alpha_{s}+(\beta_{s}+\gamma\eta_{s})\|Y\|_{\s^{\infty}_{[s,T]}}+k_s\Big)\sum_{i=1}^{n}\exp(\gamma|Y_{s}^{i}|) {\rm d}s.
\end{array}
\end{equation} 
Taking expectation conditioned on $\F_{\tau}$ for $\tau \in [t_0,t]$, we show that for $t\in [t_0,T]$,
\begin{equation}\label{eq:4.9}
\begin{array}{lll}
&& \Dis \E\big[\sum_{i=1}^{n}\exp(\gamma|Y_{t}^{i}|)|\F_{\tau}\big] \vspace{0.2cm}\\
& \leq & \Dis 2n(\exp(\gamma\|\xi\|_{\infty})+1) +\int_{t}^{T}2\gamma\big(\alpha_{s}+(\beta_{s}+\gamma\eta_{s})\|Y\|_{\s^{\infty}_{[s,T]}}+k_s\big)\E\big[\sum_{i=1}^{n}\exp(\gamma|Y_{s}^{i}|)|\F_{\tau}\big] {\rm d}s.
\end{array}
\end{equation}  
Using Gronwall's inequality, we get 
\begin{equation}\label{eq:4.10}
\begin{array}{lll}
&&\Dis \E\big[\sum_{i=1}^{n}\exp(\gamma|Y_{t}^{i}|)|\F_{\tau}\big] \vspace{0.2cm}\\
& \leq & \Dis  2n\Big(\exp(\gamma\|\xi\|_{\infty})+1\Big)\exp\Big(\int_{t}^{T}2\gamma\big(\alpha_{s}+(\beta_{s}+\gamma\eta_{s})\|Y\|_{\s^{\infty}_{[s,T]}}+k_s\big){\rm d}s\Big), \ \ t\in [t_0,T].
\end{array}
\end{equation} 
Let $\tau=t$, we have for $t\in [t_0,T]$,
$$
\sum_{i=1}^{n}\exp(\gamma|Y_{t}^{i}|)  \leq   2n\Big(\exp(\gamma\|\xi\|_{\infty})+1\Big)\exp\Big(\int_{t}^{T}2\gamma\big(\alpha_{s}+(\beta_{s}+\gamma\eta_{s})\|Y\|_{\s^{\infty}_{[s,T]}}+k_s\big){\rm d}s\Big).
$$
Using Jensen's inequality, we obtain that 
$$
\sum_{i=1}^{n}\exp(\gamma|Y_{t}^{i}|) \geq n\exp(\frac{1}{n}\sum_{i=1}^{n}\gamma|Y_{t}^{i}|) \geq n\exp(\frac{\gamma|Y_{t}|}{n}).
$$
Hence we have
\begin{equation}\label{eq:4.11}
|Y_{t}|  \leq  \Dis \frac{n}{\gamma}\log\Big(2\exp(\gamma\|\xi\|_{\infty})+2\Big)+ \int_{t}^{T}2n\Big(\alpha_{s}+(\beta_{s}+\gamma\eta_{s})\|Y\|_{\s^{\infty}_{[s,T]}}+k_s\Big){\rm d}s, \ \ t\in [t_0,T].
\end{equation} 
Since
$$
\eta_{s} \leq \eta_{s}\log(1+\eta_{s})+2,
$$
we have 
$$
\int_{0}^{T}\eta_{s}{\rm d}s \leq C_2 +2T.
$$
From the definition of $k_s$, we have
$$
\int_{0}^{T}k_{s}{\rm d}s \leq (1+\frac{2n}{\gamma})(C_2+2T)+C_2.
$$
Let
$$
C_3:=\frac{n}{\gamma}\log\big(2\exp(\gamma C_1)+2\big)+2n(1+\frac{2n}{\gamma})(C_2+2T)+4nC_2.
$$
From (\ref{eq:4.11}), we have
\begin{equation}\label{eq:4.12}
\|Y\|_{\s^{\infty}_{[t,T]}}  \leq  \Dis C_3+ \int_{t}^{T}2n(\beta_{s}+\gamma\eta_{s})\|Y\|_{\s^{\infty}_{[s,T]}}{\rm d}s, \ \ t\in [t_0,T].
\end{equation} 
Using Gronwall's inequality, we get 
\begin{equation}\label{eq:4.13}
\|Y\|_{\s^{\infty}_{[t_0,T]}}  \leq  \Dis C_3\exp\Big(\int_{t_0}^{T}2n(\beta_{s}+\gamma\eta_{s}){\rm d}s\Big) \leq C_3\exp\Big(2nC_2(\gamma+1)+4n\gamma T\Big).
\end{equation} 
Let
$$
\lambda :=C_3\exp\Big(2nC_2(\gamma+1)+4n\gamma T\Big).
$$
We get \eqref{eq:4.2}. The proof is complete.\vspace{0.2cm}

{\bf Proof of Theorem~\ref{thm:2.2}}\ \ For the number $\lambda$ given in Lemma~\ref{lem:4.2}, we have 
$$
\|\xi\|_{\infty} = \|Y\|_{\s^{\infty}_{[T,T]}} \leq \lambda.
$$
From Theorem~\ref{thm:2.1}, there exists $t_{\lambda} > 0$ which depends on constants $(n,\gamma,\delta,T,K,C_0,\lambda)$, such that BSDE \eqref{eq:1.1} has a local solution $(Y,Z) \in \s^{\infty}(\R^n) \times {\rm BMO}(\R^{n\times d}) $ on $[T-t_{\lambda},T]$. From Lemma~\ref{lem:4.2}, we obtain that 
$$
\|Y_{T-t_{\lambda}}\|_{\infty} \leq \|Y\|_{\s^{\infty}_{[T-t_{\lambda},T]}} \leq \lambda.
$$
Taking $T-t_{\lambda}$ as the terminal time and $Y_{T-t_{\lambda}}$ as the terminal value, from Theorem~\ref{thm:2.1} we know BSDE \eqref{eq:1.1} has a local solution $(Y,Z) \in \s^{\infty}(\R^n) \times {\rm BMO}(\R^{n\times d}) $ on $[T-2t_{\lambda},T-t_{\lambda}]$. Stitching the solutions we have a solution $(Y,Z) \in \s^{\infty}(\R^n) \times {\hcal^2}(\R^{n\times d})$ on $[T-2t_{\lambda},T]$ and $\|Y_{T-2t_{\lambda}}\|_{\infty} \leq \lambda$. Repeating the preceding process, we can extend the pair $(Y,Z)$ to the whole interval $[0,T]$ within finite steps such that $Y$ is uniformly bounded by $\lambda$ and $Z \in {\hcal^2}(\R^{n\times d})$. We now show that  $Z \in {\rm BMO}(\R^{n\times d})$.
Identical to the proof of inequality \eqref{eq:4.5} and \eqref{eq:4.6}, we have
\begin{equation}\label{eq:4.14}
 \Dis \eta_{s}\log(|Z_{s}|+1) \leq \frac{\gamma}{4n}\exp(-\gamma\|Y\|_{\s^{\infty}_{[s,T]}})|Z_{s}|^2+\eta_{s}(1+\gamma\|Y\|_{\s^{\infty}_{[s,T]}}+\log(\eta_{s}+1)+\frac{4n}{\gamma}),
\end{equation}
and
\begin{equation}\label{eq:4.15}
\begin{array}{lll}
\Dis u(Y_{t}^{i}) & \leq & \Dis  u(\xi^{i})-\int_{t}^{T}u'(Y_{s}^{i})Z_{s}^{i}{\rm d}W_s+\int_{t}^{T}\big(-\frac{1}{2}|Z_{s}^{i}|^{2}+\frac{1}{4n}|Z_{s}|^{2}\big) {\rm d}s \vspace{0.2cm}\\
& & \Dis +\int_{t}^{T}\gamma^{-1}\exp(\gamma|Y_{s}^{i}|)\big(\alpha_{s}+(\beta_{s}+\gamma\eta_{s})\|Y\|_{\s^{\infty}_{[s,T]}}+\hat{k}_s\big) {\rm d}s,
\end{array}
\end{equation}
where 
$$
\hat{k}_s=\eta_{s}(1+\log(\eta_{s}+1)+\frac{4n}{\gamma}).
$$
Summing $i$ from $1$ to $n$ and taking expectation conditioned on $\F_{t}$, we have
\begin{equation*}
\begin{array}{lll}
\Dis & & \Dis \frac{1}{4}\E\big[\int_{t}^{T}|Z_{s}|^{2}{\rm d}s|\F_{t}\big] \vspace{0.2cm}\\
& \leq & \Dis \gamma^{-2}n\exp(\gamma C_1) + \gamma^{-1}n\exp(\gamma \lambda)\bigg((\lambda+2)C_{2}+(\lambda\gamma+1+\frac{4n}{\gamma})(C_2+2T)\bigg).
\end{array}
\end{equation*}
Hence $Z \in {\rm BMO}(\R^{n\times d})$.

Finally, we prove the uniqueness. Let $(Y,Z)$ and $(\widetilde Y, \widetilde Z)$ be two solutions in $\s^{\infty}(\R^n) \times {\rm BMO}(\R^{n\times d})$ on $[0,T]$. Denote
$$
K_0:=\|Y\|_{\s^{\infty}_{[0,T]}}+\|\widetilde Y\|_{\s^{\infty}_{[0,T]}}+\|Z\|_{{\rm BMO}_{[0,T]}}^2+\|\widetilde Z\|_{{\rm BMO}_{[0,T]}}^2.
$$
Identical to the proof of inequality \eqref{eq:3.19}, on the interval $[T-\eps,T]$, we have
\begin{equation}\label{eq:4.16}
\begin{array}{lll} 
&&\Dis \big\|\Delta Y\big\|_{\s^{\infty}}^2+c_1^{2}\big\|\Delta Z\big\|_{{\rm BMO}(\mathbb P)}^2 \vspace{0.2cm}\\
&\leq&\Dis 64n\eps\big[\phi(2K_0)\big]^{2}\Big(T+16c_2^{2}K_0+8K_0\Big) \big\|\Delta Y\big\|_{\s^{\infty}}^{2}\vspace{0.2cm}\\
&&\Dis+96\sqrt{3}n\eps^{1-\delta}\big[\phi(2K_0)\big]^{2}\Big(1+nL_4^{2}c_2^{2}\Big)\Big(T^{\delta}+8+12L_4^{2}c_2^{2}K_0+4K_0\Big)\big\|\Delta Z\big\|_{{\rm BMO}(\mathbb P)}^{2}, 
\end{array}
\end{equation}
where $c_1$ and $c_2$ only depend on $(T,K_0)$ and the function $\phi(\cdot)$. When $\eps$ is sufficiently small, we get $Y=\widetilde Y$ and $Z=\widetilde Z$ on $[T-\eps,T]$. Repeating the preceding process within finite steps, we get the uniqueness on $[0,T]$. The proof is complete.  
\bigskip
%%%%%%%%%%%%%%%%%%%%%%%%%%%%%%%%%%%%%%%%%%%%%%%%%%%%%%%%%%%%%%%%%%%%%%%%%%%%

\end{document}